\documentclass[12pt]{article}
\usepackage{ifpdf}
\usepackage{subfigure}
\usepackage[english]{babel}
\usepackage{graphicx}
\usepackage{mathrsfs}
\usepackage{amsfonts}
\usepackage{amsmath}

\usepackage{amssymb}
\usepackage{amsmath}
\usepackage{amsthm}
\usepackage[utf8]{inputenc}
\usepackage[T1]{fontenc}
\usepackage{lmodern, tabu}
\usepackage{ae}

\ifx\pdfoutput\undefined
\usepackage[hypertex]{hyperref}
\else
\usepackage[pdftex,hypertexnames=false]{hyperref}
\fi

\newtheorem{theorem}{Theorem}[section]

\newtheorem{proposition}[theorem]{Proposition}
\newtheorem{remark}[theorem]{Remark}

\newcommand{\dys}{\frac{d \, y}{d \, s}}

\begin{document}
\thispagestyle{empty}
\vfill
\vspace{0.5cm}
\begin{center} {\Large Algebraic-based nonstandard time-stepping schemes} \end{center}
\vspace{0.5cm}
\begin{center}{\sc Lo\"ic MICHEL} \end{center}
\vspace{0.5cm}

\begin{quotation}                
\begin{center} {\bf Abstract} \end{center}
\noindent
In this preliminary work\footnote{This work is distributed under CC license \url{http://creativecommons.org/licenses/by-nc-sa/4.0/}. Email of the corresponding author: loic.michel54@gmail.com}, 
we present nonstandard time-stepping strategies to solve differential equations based on the algebraic estimation method applied to the estimation of time-derivative, which provides interesting properties of "internal" filtering. 
We consider firstly a classical finite difference method, like the explicit Euler method for which we study the possibility of using the algebraic estimation of derivatives instead of the usual finite 
difference to compute the numerical derivation. Then, we investigate how to use the algebraic estimation of derivatives in order to improve the slope predictions in RK-based schemes.
\end{quotation}

\tableofcontents

\newpage
\setcounter{page}{1}
\section{Introduction}
Ordinary differential equations and stiff differential equations \cite{Hairer, Butcher} have been studied extensively and successful methods have been proposed (e.g. \cite{Guzel, Abelman,
Enright, Jannelli, Hojjati, Darvishia, Cash, Ketcheson}), including the solver routines 'ode45' and 'ode23s' \cite{Byrne, Shampine, Ashino}.
Recent alternative methods propose for instance to use specific series, like power series \cite{ben} or Haar wavelets \cite{Hsiao, Lupik} to describe the solution of ODEs; the purpose is to substitute
the standard finite difference by the numerical properties of the series.
In this paper, we propose to extend the finite difference technique by an estimation of the derivatives using the algebraic estimation approach. 
Introduced in \cite{fliess}, the algebraic estimation method \cite{sira, fliess2} has been widely applied to many different problems of estimation that occur in dynamical 
systems. 
Some of these applications aim e.g. to reconstruct the states of dynamical systems \cite{reger, sira}, 
and a computational toolbox has been released to help processing efficiently the algebraic estimation of derivatives for particular problems \cite{Zehetner}. The algebraic estimation of derivatives
could be seen as similar to the differentiation by integration technique, investigated especially by \cite{Huang} and \cite{Wang}, and probabilistic Kriging-based method \cite{vasquez}.

Our proposed method belongs {\it a priori} to the class of nonstandard finite-difference (NSFD) methods \cite{Mickens_1, Mickens_2, Manning, Anguelov, Abelmana}, described as "powerful
numerical methods that preserve significant properties of exact solutions of the corresponding differential equations" (an interesting survey can be found in \cite{Patidar});
explicit rules to "design" NSFD schemes have been proposed in \cite{Mickens_1} \cite{Mickens_2}. Among the derivations that have been proposed (e.g. \cite{Erdogana, Ing, Anguelov_2, bru}), 
one emphasizes the contribution of \cite{bru} that describes a nonstandard finite-difference scheme for fractional systems, that uses a discrete version of the Caputo fractional derivation.

In this paper, we attempt to design simple NSFD schemes, that are derived from the algebraic estimation technique  in order to evaluate the discrete derivatives of first order ODEs considering
an "Euler framework" and a "RK framework". In an "Euler-like scheme", the algebraic estimation of the derivatives is mainly used to build a multi-step scheme, where the evaluation of the derivative depends 
also on the past values of the solution. In a RK-like scheme, the estimation of the local slope, that is usually
"measured" as an average of several local slopes, is performed using the properties of filtering of the algebraic technique.

The paper is structured as follow. Section 2 described the proposed methods, including a brief review of the algebraic estimation method. Some concluding remarks can be found in Section 3. 

\newpage
\section{Outline of the method}

Consider an ordinary differential equation (ODE) such as:
\begin{equation}\label{eq:ODE_fond} 
\frac{d \, y(t)}{d \, t} = f(y(t), u(t)), \qquad t \in [ 0, \, t_f], \quad y(0) = y_0.
\end{equation}
\noindent
The quantities $u$ and $y$ represent respectively the input and the solution of (\ref{eq:ODE_fond}). The corresponding usual discrete explicit Euler scheme reads:

\begin{equation}\label{eq:EDP_discrete} 
\left. \frac{d \, y(t)}{d \, t} \right|_{k+1} =  \frac{y_{k+1} - y_k}{h} \approx f(y_k, u_k), \qquad k \in \mathbb{N}, \quad y(0) = y_0
\end{equation}

\noindent
where $k$ is the sampled time, $h$ is the time-step, and $y_k$, $y_{k+1}$ are respectively the solution of (\ref{eq:EDP_discrete}) at the discrete instants $t_k$ and $t_{k+1}$. The sampled are supposed equally
distributed i.e. $t_k - t_{k-1}$ is constant for all $k$.

\vspace{0.5cm}
We identified two major strategies to build time-stepping schemes, based on the rules that help designing NSFD schemes and the algebraic estimation technique. The "Euler-like"(or "multi-step") strategy follows the principle of
the Euler method, and the "RK-like" strategy is based on the general Runge-Kutta method \cite{Shampine}.

\subsubsection*{Algebraic estimator of derivatives}

Using the previous notations, consider a function $g(x)$ and define the algebraic estimator of derivative $\mathcal{D}(g)_n$ of order $n$ such as:

\begin{equation}\label{eq:estimator}
\begin{array}{c}
\displaystyle{\left. \frac{d \, g(x)}{d \, x} \right|_{k+1}  \approx \mathcal{D}(g)_n = \frac{\alpha_0 g_{k+1} + \alpha_1 g_k + \cdots + \alpha_n g_{k-n}}{\phi(h)}}, \quad  \\
\\[0.01cm]
\phi(h) = h + O(h^2) \,\, \mathrm{ as } \,\,  h \rightarrow 0 \quad \mathrm{where} \quad h = x_{k} - x_{k-1} \quad \mathrm{is} \,\, \mathrm{constant}
\end{array}
\end{equation}

\noindent
The coefficients $\alpha_0, \, \alpha_1, \cdots, \alpha_n$ are real coefficients. 
These coefficients and the function $\phi(h)$ are defined by (\ref{eq:coeff_def}) in the proposition (\ref{prop_1}) as calculating steps of the proposed time-stepping 
schemes\footnote{The algebraic estimator of derivative is defined with the highest degree of $g$ that is equal to $k+1$. Depending on the application, this degree can be decreased to $k$.}. We have:

\begin{equation}\label{eq:coeff_def}
\begin{array}{c}
 \alpha_0 = T \\
 \alpha_i = 2(T - 2 h i), \, i \in \{ 1 \cdots \eta-1 \} \\
 \alpha_n = (T - 2 h \eta) \\
 \displaystyle{\frac{1}{\phi(h)} = \frac{3Kh}{T^3} = \frac{3Kh}{(\eta h)^3}}
 \end{array}
\end{equation}

\newpage
\subsection{Euler-like Algebraic-NSFD scheme}

\subsubsection{Definition}

The proposed "Euler-like" Algebraic-NSFD scheme aims to extend the finite difference in (\ref{eq:ODE_fond}) by the algebraic estimator $\mathcal{D}(f)_n$.
First, one proposes a "symbolic" nonstandard scheme, that is of the form:

\begin{equation}\label{eq:gen_scheme} 
\begin{array}{c}
\mathcal{D}(f)_n \approx f(y_k, u_k)  \\
\\[0.01cm]
\mathrm{i.e.} \quad \displaystyle{\frac{\alpha_0 y_{k+1} + \alpha_1 y_k + \cdots + \alpha_n y_{k-n+1}}{\phi(h)}} \approx f(y_k, u_k), \, k \in \mathbb{N^{*+}}, \quad y_0 = y(0).
\end{array}
\end{equation}

\noindent
A forward scheme can be easily deduced from the proposed general scheme (\ref{eq:gen_scheme}):

\begin{equation}\label{eq:gen_scheme_fond} 
y_{k+1} =  - \frac{ \alpha_1 y_k + \cdots + \alpha_n y_{k-n+1}}{\alpha_0} + \frac{\phi(h)}{\alpha_0} f(y_k, u_k) , \qquad k \in \mathbb{N}^{*+}, \quad y_0 = y(0).
\end{equation}
\noindent
As presented in \cite{bru} regarding the Caputo fractional derivative, due to the nonlocal nature of the algebraic-based derivative operator, the discrete representation of the derivative must take into account a part of the past history of the 
solution. The number $n$ of the involved sampled solutions $y_{k}, \cdots, y_{k-n}$ defines a window that characterizes the "precision" of the derivative estimation\footnote{This estimation window implies 
that the $n$ initial sampled solutions are obviously not known at the beginning of the algorithm. To initialize the estimation window, one may consider e.g. using the classical finite difference scheme.}.

\vspace{0.5cm}
Then, in the following proposition, we formalize the proposed Euler-like nonstandard time-stepping scheme based on the algebraic estimation framework.
 
\begin{proposition}
\label{prop_1}
Consider the following nonstandard numerical scheme associated to the ODE (\ref{eq:ODE_fond}), that verifies the general scheme (\ref{eq:gen_scheme}):
\begin{equation}\label{eq:EDP_NSFD} 
\begin{array}{c}
\displaystyle{ K \frac{3h}{T^3} \left\{ \sigma_0 + \sum_{j=1}^{\eta-1} 2y_{k-\eta+j+1}(T - 2 j h) \right\} } \approx f(y_k, u_k), \qquad k \in \mathbb{N}^{*+}, \quad y_0 = y(0) \\
\\[0.01cm]
\hbox{ \rm{with} } \displaystyle{T = \eta h > 0; \, \eta \in \mathbb{N}^{*+} } \hbox{ \rm{and} }  \sigma_0 = y_{k-\eta+1} T + y_{k+1} (T - 2 \eta h)
\end{array}
\end{equation}
\noindent
where $k$ is the sampled time, $h$ is the time-step, $K$ is a real constant, and $T$ is a multiple of $h$ that
characterizes the "low filtering" property of the algebraic derivative (see $\S$1 in \ref{first_para}). This scheme is called {\it Euler-Algebraic-NSFD} scheme, or simply {\it E-A-NSFD} scheme with 
a window $T$. 

\end{proposition}

\begin{proof}
\underline{Hypothesis} We consider solving the ODE (\ref{eq:ODE_fond}), for which the solution $y(t)$ is assumed, in the time domain, to be locally represented by a linear function of the time i.e.:

\begin{equation}\label{eq:local_est}
y(t) = a_0 + a_1 t
\end{equation}

\noindent
where, in particular, the coefficient $a_1$ is calculated from the algebraic estimation technique. The lowest degree of "time-filtering" is considered.

\vspace{0.5cm}
\underline{Step 1 - Algebraic Derivation} The general technique to perform the estimation of derivatives using the algebraic estimation strategy is mainly described in \cite{fliess, sira}. It
allows estimating the coefficient $a_1$ in (\ref{eq:local_est}). We consider only the first algebraic derivative of $f$ according to the definition of the initial-value problem (\ref{eq:ODE_fond}). 
Transform first (\ref{eq:local_est}) in the Laplace domain, then take the derivative $d / ds$:

\begin{equation}
 y(s) = \frac{a_0}{s} + \frac{a_1}{s^2} \Longleftrightarrow s y(s) = a_0 + \frac{a_1}{s}
\end{equation}

\begin{equation}\label{eq:alg_deriv}
y(s) + s \dys  = - \frac{a_1}{s^2}
\end{equation}

\vspace{0.5cm}
\underline{Step 2 - Back to the time domain} Multiply first (\ref{eq:alg_deriv}) by $s^{-p}$:

\begin{equation}\label{eq:alg_deriv_2}
s^{-p} y(s) + s^{-p+1} \dys = - s^{-p-2} a_1
\end{equation}

\noindent
then using the Cauchy formula\footnote{Some elements of proof can be found in \cite{delpoux} (pp.22-23).} applied to each term of (\ref{eq:alg_deriv_2}):

\begin{equation}\label{eq:cauchy}
\int^{(\beta)} \frac{1}{s^{\alpha}} \frac{d^{\beta} x(\tau)}{d \, \tau^{\beta}} d^{\beta} \, \tau = \frac{1}{(\alpha - 1) !} \int_{0}^{t} (t- \tau)^{\alpha-1} (-1)^{\beta} \tau^{\beta} x(\tau) d \tau
\end{equation}

\noindent
we deduce\footnote{The expression of $a_1$ corresponds to the direct application of (\ref{eq:cauchy}), whose integral is evaluated over a moving windows $T$.} the expression of $a_1$ that corresponds 
to the derivative of $y$ through (\ref{eq:local_est}):

\begin{equation}\label{eq:deriv_gen_exp}
a_1 = \frac{ \displaystyle{\frac{1}{(p - 1)!}  \int_{t-T}^{t} (t - \tau)^{p-1} y( \tau) d \, \tau - \frac{1}{(p-2)!} \int_{t-T}^{t} (t - \tau)^{p-2} \tau x( \tau) d \, \tau}}{\displaystyle{\frac{1}{(p +1)!} \int_{t-T}^{t}(t - \tau)^{p+1} d \, \tau}}
\end{equation}

\noindent
Since we consider the lowest degree of "time-filtering", we take $p = 2$. 

\newpage
\noindent
Therefore, (\ref{eq:deriv_gen_exp}) reads:

\begin{equation}\label{eq:_3}
a_1 = - \frac{6}{T^3} \displaystyle{ \int_{t-T}^{t} (t - 2 \tau) y( \tau) d \, \tau}
\end{equation}

\vspace{0.5cm}
\underline{Step 3 - Integral expression of the numerical scheme} Replacing in (\ref{eq:ODE_fond}) the time-derivative operator by (\ref{eq:_3}) gives an integral expression of the proposed E-A-NSFD scheme:

\begin{equation}\label{eq:_4}
-\frac{6}{T^3} \displaystyle{ \int_{t-T}^{t} (t - 2 \tau) y( \tau) d \, \tau} \approx f(y(t), u(t))
\end{equation}

\vspace{0.5cm}
\underline{Step 4 - Finally, the discretized version...} The usual Trapezoidal scheme allows integrating numerically and thus providing a discretized expression of (\ref{eq:_4}). We have:

\begin{equation}\label{eq:discrete_deriv}
a_1 = - K \frac{3h}{T^3} \left\{ y_{k-\eta+1} T + y_{k+1} (T - 2 \eta h)  + \sum_{j=1}^{\eta-1} 2y_{k-\eta+j+1}(T - 2 j h) \right\}
\end{equation}

\noindent
which, as a result, gives the proposed NSFD scheme (\ref{eq:EDP_NSFD}) when substituted in (\ref{eq:EDP_discrete}) according to the rule \#2 of the design methodology of NSFD schemes\footnote{Since 
the complete finite difference scheme in (\ref{eq:EDP_discrete}) is {\it a priori} replaced by a weighted sum of the discrete values $y_k$, we deduce that the 
rule \#2, regarding the design of NSFD scheme \cite{Patidar}, is not completely verified (since only the $"h"$ term has to be changed). A future work should investigate this ambiguity. 
}
\cite{Patidar}. 
Depending on the value of $T$, the total gain of $|a_1|$ could be different, that is why 
we multiply $a_1$ by a positive constant $K$ that should be estimated at the initial instants\footnote{In simulation, it is observed that $K \rightarrow 1$ when $T$ increases.}. 
Note that if $\eta = 1$, then $T=h$ which cancel the summation term $\sum_{j=1}^{\eta-1} 2y_{k-\eta+j+1}(T - 2 j h)$.

\vspace{0.5cm}
\underline{Step 5 - In addition, identification with the general scheme} 
The coefficients $\alpha_i, \, i \in \{1 \cdots n \}$ can be identified to the general form (\ref{eq:gen_scheme}). We have:

\begin{equation}\label{eq:coeff_def}
\begin{array}{c}
 \alpha_0 = T \\
 \alpha_i = 2(T - 2 h i), \, i \in \{ 1 \cdots \eta-1 \} \\
 \alpha_n = (T - 2 h \eta) \\
 \displaystyle{\frac{1}{\phi(h)} = \frac{3Kh}{T^3} = \frac{3Kh}{(\eta h)^3}}
 \end{array}
\end{equation}
\end{proof}

As an example of the Prop. \ref{prop_1}, the simplest version of the proposed E-A-NSFD is a particular case where $T = h$ that verifies (\ref{eq:EDP_discrete}). 

\newpage
\noindent
To illustrate how the E-A-NSFD scheme is written, consider the simple ODE:
 
\begin{equation}\label{eq:ODE_ex} 
\frac{d \, y(t)}{d \, t} = 5 y(t) + u(t), \qquad t \in [ 0, \, t_f], \quad y(0) = y_0
\end{equation}

\begin{itemize}

\item Considering $\eta = 3; h = 10^{-5}$ ($\Psi_3$ is a real constant), the E-A-NSFD scheme of (\ref{eq:ODE_ex}) reads: 
\begin{equation}
\Psi_3 (3u_{k-2}+2y_{k-1}-2y_k-3y_{k+1}) \approx 5 y(t) + 50 u(t)
\end{equation}

\item Considering $\eta = 5; h = 10^{-5}$ ($\Psi_5$ is a real constant), the E-A-NSFD scheme of (\ref{eq:ODE_ex}) reads: 
\begin{equation}
 \Psi_5 (5y_{k-4}+6y_{k-3}+2y_{k-2}-2y_{k-1}-6y_k-5y_{k+1}) \approx 5 y(t) + 50 u(t)
\end{equation}
\end{itemize}

\noindent
As stated in Prop. \ref{prop_2}, the distribution of the signs in these two cases is equally distributed. 

\subsubsection{Some properties}
\label{first_para}

This strategy could been seen as a receding horizon approach because the interval of integration $T$ corresponds to the window from which the estimation of the derivative is performed. Therefore,
the constant $T$ should  be chosen small enough in order to evaluate the  estimate  within  an  acceptable  short  delay, but large enough, in order to ensure a good low pass filtering \cite{Zehetner}.

Since the whole set of coefficients $\alpha_i$ does not depend of the discrete solution $y_k$, they need only to be evaluated at the initialization step of the scheme.

The following simple proposition gives an interesting property of the general scheme (\ref{eq:gen_scheme}) concerning the distribution of the signs of the set of $\alpha_i$.

\begin{proposition}
\label{prop_2}
The distribution of the signs (over the whole set of the coefficients $\alpha_i, \, \, i = 0 \cdots \eta$) verifies: 
$\mathrm{sign}(\alpha_0, \cdots, \alpha_{\lfloor \eta/2-1 \rfloor}) = -\mathrm{sign}(\alpha_{\lfloor \eta/2 \rfloor}, \cdots, \alpha_{\eta})$.
\end{proposition}

\begin{proof}
From Prop. \ref{prop_1}, we have $\alpha_i = \gamma (T - 2 h i), \, i = 0 \cdots \eta, \, \gamma \in \mathbb{N}^{*+}$. To get the number of $\alpha_i$ that are e.g. positive, 
check the inequality $(T - 2 h i) > 0$. This gives $i < \lfloor\frac{\eta}{2} \rfloor$.
\end{proof}

The E-A-NSFD scheme may become unstable for high $\eta$ \footnote{The window $T$ of the algebraic estimation 
may be seen as an "internal" time-delay inside the ODE (\ref{eq:EDP_discrete}) that may "deform" the value of the time-derivative estimation. Further works may investigate the conditions of stability according 
to the length of the window $T$ and the coefficient $K$.}.

\subsubsection{Toward possible generalizations}
\label{generaliz}
One may consider a possible generalization of the scheme (\ref{eq:EDP_NSFD}) using a more general expression $\phi(h)$ instead of the constant time-step $h$. The scheme can be rewritten:
 
\begin{equation}
\frac{3 K \phi}{T^3} \left\{ \sigma_0 + \sum_{j=1}^{\eta-1} 2y_{k-\eta+j+1}(T - 2 j \phi) \right\}\approx f(y(t), u(t))
\end{equation}
 
\noindent
with $\phi(h) = h + O(h^2)$, as $h \rightarrow 0$.
 
 
\subsection{RK-like Algebraic-NSFD scheme}

The family of RK methods is described generally as an Euler scheme for which multiple estimations of the local slopes around $f(y_k)$ are computed in such manner that the solution increment $y_{k+1}$ is based 
upon a weighted average of these multiple estimations. The general RK scheme reads:

\begin{equation}\label{eq:RK_gen}
 y_{k+1} = y_k + h \sum_{i = 1}^s b_i q_i
\end{equation}

\noindent
where $q_1, q_2, \cdots, q_s$ are the estimations of the different slopes through $f$. The number $s$ characterizes the choice of the method in the RK family.

\subsubsection{Definition}

The proposed "RK-like" Algebraic-NSFD scheme aims to substitute the weighted sum $\sum_{i = 1}^s b_i q_i$ in (\ref{eq:RK_gen}) by the algebraic estimator $\mathcal{D}(f)_n$.

First, one proposes a "symbolic" nonstandard scheme, that is of the form:

\begin{equation}\label{eq:RK_gen_scheme} 
\begin{array}{c}
y_{k+1} = y_k + \phi(h) f(y_k + \phi(h) [f(y_k) + \phi(h) \mathcal{D}(f)_n]) \quad \mathrm{i.e.} \\
\displaystyle{\frac{y_{k+1} - y_k}{\phi(h)} = f(y_k + \phi(h) [f(y_k) + \phi(h)(\alpha_0 f(y_k + \delta_{0}) + \cdots + \alpha_n f(y_{k} + \delta_{n}))])} \\
 k \in \mathbb{N^{*+}}, \quad y_0 = y(0) 
\end{array}
\end{equation}

\noindent
where the coefficients $\delta_{0}, \delta_{1}, \cdots, \delta_{n}$ are real coefficients such as for all $i \in [0, \,  n], \delta_i \in [\delta_{min}, \, \delta_{max}]$; 
$\alpha_0 \cdots \alpha_n$ are real coefficients and $\phi(h) = h + O(h^2)$, as $h \rightarrow 0$.

\vspace{0.5cm}
Then, in the following proposition, we formalize the proposed RK-like nonstandard time-stepping scheme based on the algebraic estimation framework.
 
\begin{proposition}
\label{prop_3}
Consider the following nonstandard numerical scheme associated to the ODE (\ref{eq:ODE_fond}), that verifies the general scheme (\ref{eq:RK_gen_scheme}):
\begin{equation}\label{eq:EDP_NSFD_RK} 
\begin{array}{c}
\displaystyle{\frac{y_{k+1} - y_k}{\phi(h)} = f\left(y_k + \frac{h}{2} \left[f(y_k) + \frac{3Kh^2}{T^3} \left( \sigma_0 + \sum_{j=1}^{\eta-1} 2f(y_{k} + \delta_{j})(T - 2 j h) \right) \right] \right) }, \\
\\[0.01cm]
\hbox{ \rm{with} } \displaystyle{T = \eta h > 0; \, \eta \in \mathbb{N}^{*+} } \hbox{ \rm{and} }  \sigma_0 = f(y_{k} + \delta_{0}) T + f(y_{k} + \delta_{\eta}) (T - 2 \eta h), \\
\\[0.01cm]
k \in \mathbb{N^{*+}}, \quad y_0 = y(0) 
\end{array}
\end{equation}
\noindent
where $k$ is the sampled time, $h$ is the time-step, $K$ is a real constant, and $T$ is a multiple of $h$ that
characterizes the "low filtering" property of the algebraic derivative (see $\S$1 in \ref{first_para}). 
This scheme is called {\it RK-Algebraic-NSFD} scheme, or simply {\it RK-A-NSFD} scheme with a window $T$. 
\end{proposition}

\begin{proof}
\underline{Hypothesis} We consider solving the ODE (\ref{eq:ODE_fond}), for which the function $f$ is assumed, to be locally represented by a linear function i.e.:

\begin{equation}\label{eq:local_est_}
f(x) = a_0 + a_1 x
\end{equation}

\noindent
where, in particular, the coefficient $a_1$ is calculated from the algebraic estimation technique. The lowest degree of "time-filtering" is considered. The purpose is to evaluate the 
local value of $a_1$ for $x \in [x_{k-1}, \, x_k]$, for all $k \geq 1$.

\vspace{0.5cm}
\underline{Step 1 - From the discretized version of the derivative estimation} From Prop. \ref{prop_1}, the expression of $a_1$ in the discrete domain is given by:

\begin{equation}\label{eq:a_1__RK}
a_1 = - K \frac{3h}{T^3} \left\{ f(x_{k-\eta}) T + f(x_{k}) (T - 2 \eta h)  + \sum_{j=1}^{\eta-1} 2f(x_{k-\eta+j})(T - 2 j h) \right\}
\end{equation}
\noindent
The set of values $x_{k-\eta}, x_{k-\eta+1}, \cdots, x_{k}$ is, according to (\ref{eq:estimator}), theoretically a regular grid, but since the interval $[x_{k-1}, \, x_k], \, k \geq 1$ is very small, we
can consider evaluating $\eta+1$ values of $f(y_k + \delta_j)$ where for all $i \in [0, \,  \eta+1], \delta_i \in [\delta_{min}, \, \delta_{max}]$. In particular, the $\delta_i$ coefficients can 
be random numbers. We have:

\begin{equation}\label{eq:a_1_RK}
a_1 = - K \frac{3h}{T^3} \left\{ f(x_k+ \delta_0) T + f(x_{k} + \delta_{\eta}) (T - 2 \eta h)  + \sum_{j=1}^{\eta-1} 2f(x_k+\delta_j)(T - 2 j h) \right\}.
\end{equation}

\newpage
\underline{Step 2 - Definition of an RK2 scheme with algebraic estimation of $f$} To build the proposed NSFD scheme, we start from the standard RK2 scheme. 

\noindent
Considering the evolution of the solution $y$ between the steps $y_{k}$ and $y_{k+1}$, for any $k \geq 1$. 
\noindent
The first step of the RK2 scheme, called {\it the prediction}, is to evaluate the solution $y_{k+1/2}$ in the middle of the step i.e.:

\begin{equation}\label{eq:RK2_1}
y_{k+1/2} = y_k + \frac{h}{2} f(y_{k})
\end{equation}

\noindent
that gives, via (\ref{eq:EDP_discrete}), the derivative of $y$ at the middle of the time-step:

\begin{equation}\label{eq:RK2_2}
\left. \frac{d \, y}{d \, t} \right|_{k+1/2} = f(y_{k+1/2}) = f\left(y_k + \frac{h}{2} f(y_{k}) \right).
\end{equation}

\noindent
Then, the second step, called {\it the correction}, is to evaluate the solution at the end of the step $y_{k}$ taking into account the derivative $\left. \frac{d \, y}{d \, t} \right|_{k+1/2}$ i.e.:

\begin{equation}\label{eq:RK2_3}
y_{k+1} = y_k + h f(y_{k+1/2}).
\end{equation}

\noindent
Instead of considering a single evaluation of the derivative at the middle of the time-step (or multiple evaluations of the derivatives like e.g. in the RK4 case), 
the proposed modification of the RK2 strategy is to estimate multiple local values of the slopes via $f$ between the steps $y_{k}$ and $y_{k+1}$ using the algebraic estimator $\mathcal{D}(f)_n$.
First, one generalizes the "prediction" step (\ref{eq:RK2_1}) that is rewritten:

\begin{equation}\label{eq:RK2_4}
y_{k+1/2} = y_k + \frac{h}{2} <f(y_k)>_h
\end{equation}

\noindent
where we denote $<f(y_k)>_h$ the averaged value of $f(y)$ between the steps $k$ and $k+1$ (in other words, over the time-step $h$, regarding the notations)\footnote{Regarding the RK4 algorithm (with no time dependence), we
have: $<f(y_k)>_h = \frac{1}{6}(f(y_k) + 2 (f(y_k + \frac{k_1}{2}) + 2 f(y_k + \frac{k_2}{2})  + f(y_k + k_3))$ where $k_1$, $k_2$ and $k_3$ are the estimated local slopes.}.

Before rewriting the "correction" step, one needs the formal definition of the derivative of a function $g(x)$ taken at the point $x=a$:

\begin{equation}\label{eq:nb_derive}
\left. \frac{d \, g}{d \, x} \right|_{x = a} = \lim_{h \rightarrow 0} \frac{g(a+h)-g(a)}{h} \approx \frac{g(a+h)-g(a)}{h} \quad \mathrm{ for } \, h \, \mathrm{very} \, \mathrm{small}
\end{equation}

\noindent
that allows defining the quantity $<f(y_k)>_h$, by substituting (\ref{eq:a_1_RK}) in (\ref{eq:nb_derive}):


\begin{equation}\label{eq:deriv_a1}
 <f(y_k)>_h \approx \left. f(y_k) + \frac{d \, f(y) }{d \, y} \right|_{k}  h = f(y_k)  + h a_1
\end{equation}

\noindent
Finally, from (\ref{eq:RK2_4}) and (\ref{eq:deriv_a1}), the correction step is thus rewritten\footnote{A simple approximation in (\ref{eq:RK2_3}) can be considered: since we assumed that
$f$ is locally described as a first order polynomial (\ref{eq:local_est_}), it could be possible to approximate $f(y_k)$ by $\gamma y_k$ where $\gamma$ is a real constant number that may be estimated at
each time-step.}:

\begin{equation}\label{eq:RK2_5}
y_{k+1} = y_k + h f(y_{k+1/2}) = y_k + h f\left(y_k + \frac{h}{2} <f(y_k)>_h \right)     
\end{equation}

\begin{equation}\label{eq:RK2_6}
\Longleftrightarrow  y_{k+1} =  y_k + h f\left(y_k + \frac{h}{2} (f(y_k) + h a_1) \right).  
\end{equation}

\noindent
As a result, (\ref{eq:RK2_6}) gives the proposed NSFD scheme (\ref{eq:EDP_NSFD_RK}) according to the rule \#3 of the design methodology of NSFD scheme.
In the same manner as reported in Prop.(\ref{prop_1}), depending on the value of $T$, the total gain of $|a_1|$ could be different, that is why we multiply $a_1$ by a positive constant 
$K$ that should be estimated at the initial instants.
\end{proof}

Since $f(y_{k+1}) - f(y_k)$ represents the slope that locally "drives" the discrete differential equation (\ref{eq:EDP_discrete}), we are looking for a good estimation of the predicted slope
through $f$ at a point $y$, that is located  between $f(y_k)$ and $f(y_{k+1})$. Considering smooth enough the function $f(x)$ where $x \in [y_{k}, \, y_{k+1}]$, the purpose of the operator 
$\mathcal{D}(f)_n$ is to estimate the derivative of $f$ in this 
particular interval in such manner that an "averaged" value of $f(x)$ for $x \in [y_{k}, \, y_{k+1}]$ can be deduced. 

\begin{remark}
 This scheme (\ref{eq:RK_gen_scheme}) is the dual form of the Euler-like scheme (\ref{eq:gen_scheme}), where the algebraic estimator is located on the left side of (\ref{eq:ODE_fond}) and is directly "connected" to the 
discrete values of $y$. In this RK-like scheme, the algebraic estimator is located on the right side of (\ref{eq:ODE_fond}) and is connected to the 
discrete values of $y$ through $f$.
\end{remark}

\begin{remark}
 The lowest degree of "time-filtering" has been considered in the proof of the Prop. \ref{prop_1} \ref{prop_3}. More investigations would clarify the impact of higher degrees on 
 the precision of the solution and the stability of the proposed schemes. Moreover, this could be an essential factor regarding the possibility of simulating dynamical models with noisy signals 
 (see Rem.  \ref{rem_stoc}).
\end{remark}

\subsubsection{Some properties and comparison with the E-A-NSFD scheme}

 The main difference in the utilization of the algebraic operator between the two proposed schemes is the "management" of the sampled $y$ solution. 
 In the Euler-like scheme, the window $T$ is a {\it moving} window that aims to performs the derivative estimation while the scheme runs; the window is thus initialized only at the beginning of the scheme. 
 In the RK-like scheme, the window $T$ is a {\it static} window, that aims to perform the derivative 
 estimation between two steps $[y_{k}, \, y_{k-1}]$; the window is thus initialized at each time-step and the number of evaluations of $f$ is proportional to the length of the window $T$ (especially, to
 compute \ref{eq:a_1_RK}).
 
 The Euler-like scheme may require {\it a priori} few samples of $y$ in order to compute $\mathcal{D}(f)_n$ while preserving {\it a priori} the global stability of the scheme. At the opposite, the RK-like scheme
 requests many samples due to the fact that a random process is involved in the evaluation of $f$ / estimation of the slopes. In such case, the average of the estimated slopes is performed naturally by the
 filtering property of the algebraic estimator.
 
 As presented in $\S$ \ref{generaliz}, the same type of generalization can be applied by substituting the time-step $h$ by a function $\phi$ that follows the rule 
 $\#$2 of the NSFD scheme design. Moreover, more sophisticated schemes can be built e.g.: association of the E-A-NSFD scheme and the RK-A-NSFD scheme; inclusion of the derivative estimator 
 inside high order RK scheme in order to refine to precision of the estimated slopes...
 
 \begin{remark}
 \label{rem_stoc}
  Taking into account the properties of filtering that provide the algebraic estimation framework, we assume that it could be possible to simulate dynamical models that involve noisy signals, like e.g. models
  in electronics that sometimes can not be dissociated from the physical noise \cite{sick}.
 \end{remark}

\section{Conclusion and further work}

We described in this paper, nonstandard finite-difference schemes that use the algebraic estimation framework in order to compute an estimate of the derivatives. We derived
two algebraic-based nonstandard schemes. The approach used for the first scheme is similar to \cite{bru}, which 
derives a nonstandard scheme from the Caputo derivative definition in order to solve fractional dynamical systems, and the second scheme aims to extend the Runge-Kutta method.
Further work include:

\begin{itemize}
 \item extensive tests of practical problems (e.g. using the set of practical problems presented in \cite{IVP});
 \item a complete stability study in order to characterize the stability domains according to the stiffness of the ODE;
 \item the application to higher order ODEs and utilization of the dedicated toolbox \cite{Zehetner} to systematize the proposed A-NSFD procedure;
 \item the application to fractional differential equations.
\end{itemize}

\section*{Acknowledgement}

The author is sincerely grateful to Dr. Edouard Thomas for his strong guidance and his valuable comments that improved this paper.

\end{document}